\documentclass[11pt]{amsart}

\usepackage{euscript}
\usepackage{amsmath}
\usepackage{amsfonts}
\usepackage{amssymb}
\usepackage{amsthm}
\usepackage{epsfig}
\usepackage{epstopdf}
\usepackage{url}
\usepackage{array}
\usepackage{amssymb}\usepackage{amscd}
\usepackage{epic}
\usepackage{eufrak}
\usepackage{tikz,underscore}
\usepackage{tikz-cd}
\usepackage{float}
\usepackage{subfigure}
\usepackage{tkz-euclide}

\usepackage{graphicx}
\usetkzobj{all}

\numberwithin{equation}{section}

\theoremstyle{plain}
\newtheorem{theorem}{Theorem}[section]
\newtheorem{lemma}[theorem]{Lemma}

\newtheorem{thm}[equation]{Theorem}
\newtheorem{cor}[equation]{Corollary}

\theoremstyle{definition}

\newtheorem{remark}[theorem]{Remark}
\newtheorem{rem}[theorem]{Remark}
\newtheorem{definition}[theorem]{Definition}

\DeclareMathOperator{\Isom}{{\mathrm Isom}}

\DeclareMathOperator{\diam}{{\mathrm diam}}

\def\N{\mathbb N}
\def\R{\mathbb R}
\def\Z{\mathbb Z}
 
\def\H{\mathbb{H}}

\def\Ga{\Gamma}

\def\La{\Lambda}
\def\Om{\Omega}
\def\al{\alpha} 
\def\eps{\epsilon}
\def\ga{\gamma}
\def\la{\lambda}

\def\geo{\partial_{\infty}}

\title{Hausdorff dimension of non-conical limit sets}

\author{Michael Kapovich}
\address{M.K.: Department of Mathematics, UC Davis, One Shields Avenue, Davis CA 95616, USA}

\email{kapovich@math.ucdavis.edu}

\author{Beibei Liu}
\address{B.L.: Max Planck Institute of Mathematics, Vivatsgasse 7, 53111 Bonn, Germany}
\email{bbliumath@gmail.com}

\date{September 17, 2019}

\begin{document}
\begin{abstract}
Geometrically infinite Kleinain groups have nonconical limit sets with the cardinality of  the continuum. In this paper, we construct  a geometrically infinite Fuchsian group such that the Hausdorff dimension of the nonconical limit set equals zero. For finitely generated, geometrically infinite Kleinian groups, we prove that the Hausdorff dimension of the nonconical limit set  is positive.
\end{abstract}

\maketitle

\section{Introduction}

Consider a discrete group $\Gamma$ of isometries of the $n$-dimensional hyperbolic space $\H^{n}$, i.e. a Kleinian group. The limit set $\Lambda(\Gamma)$ of $\Ga$ is the accumulation set on the ideal boundary $S^{n-1}=\geo \H^n$  of the $\Ga$-orbits in $\H^{n}$. The  notion of \emph{geometric finiteness} of $\Ga$ (going back to Ahlfors, \cite{Ah}) reflects the geometry of the quotient space $\H^n/\Ga$ and the dynamics of the $\Ga$-action on the limit set. For instance, see \cite{BM}, 
a Kleinian group $\Gamma$ is geometrically finite if and only if every limit point of $\Gamma$ is either a conical limit point or a bounded parabolic fixed point.  Here, a point $\xi\in \Lambda(\Gamma)$ is \emph{a conical limit point } of $\Gamma$ if  one, equivalently, every geodesic ray $\mathbb{R}_{+}\rightarrow \H^{n}$ asymptotic to $\xi$ projects to a non-proper map $\mathbb{R}_{+}\rightarrow M=\H^{n} / \Gamma$. In contrast, a point $\xi\in \Lambda(\Gamma)$ is a \emph{non-conical limit point} if for some (equivalently, every) 
geodesic ray asymptotic to $\xi$, the projection of this ray to the quotient space $M=\H^{n}/ \Gamma$ is proper, i.e.  eventually leaves any 
compact subset. There are only countably many parabolic fixed points in $\Lambda(\Gamma)$. Thus, for each geometrically finite group, 
 the Hausdorff dimension $\dim(\Lambda_{c}(\Gamma))$ of the set $\Lambda_{c}(\Gamma)$ of conical limit points equals the Hausdorff dimension  of $\Lambda(\Gamma)$, while the Hausdorff dimension of the nonconical limit set $\Lambda_{nc}(\Gamma)= \La(\Ga) \setminus \La_{c}(\Ga)$ is $0$. Note that it was proven by  Bishop and Jones \cite{Bi} that $\dim(\Lambda_{c}(\Gamma))$ equals the {\em critical exponent} 
 $\delta(\Gamma)$ for any non-elementary Kleinian group. 

It was proven by Bishop in \cite{Bi2} that for each \emph{geometrically infinite} torsion-free 
Kleinian group $\Gamma< \Isom(\H^3)$, the set of nonconical limit points $\La_{nc}(\Ga)$  has the cardinality of the continuum \cite{Bi2}; this result was generalized to Kleinian subgroups of $\Isom(\H^n)$ (and, more generally, of isometry groups of Hadamard manifolds of pinched negative curvature)  in our prior work \cite{KL}.  This, however, leaves open the question of the {\em measure-theoretic size} of the nonconical limit set. For discrete isometry groups of the  hyperbolic plane, Fern\'andez and Meli\'an proved the following.

\begin{theorem}\cite{FM}
\label{FMCH}
If $\Gamma<\Isom(\H^{2})$ is a torsion-free, geometrically infinite Fuchsian group such that the limit set is $S^{1}$, then  
$$\dim(\Lambda_{nc}(\Gamma))=\dim(\Lambda(\Gamma))=1.$$
\end{theorem}

\begin{remark}
Meli\'an, Rodr\'iguez and Tour\'is later generalized this result to discrete isometry groups of simply-connected complete Riemannian surfaces $X$ with curvature satisfying $-k^{2}\leq K_X\leq -1$,  proving that the visual dimension of the nonconical limit set is greater than or equal to $1$ if the limit set is the entire ideal boundary $\geo X$, 
\cite{MRT}. 
\end{remark}

A Kleinian group $\Gamma<\Isom(\H^{n})$  is said to be of the \emph{first kind} if $\Lambda(\Gamma)=S^{n-1}$. 
Otherwise, $\Ga$ is of the \emph{second kind}.  Theorem \ref{FMCH}, thus, identifies the Hausdorff dimension of the nonconical limit set of torsion-free, geometrically infinite Fuchsian groups of the first kind. For finitely generated  geometrically infinite Kleinian groups $\Ga< \Isom(\H^3)$ of the second kind, Bishop and Jones proved:
 
 \begin{theorem}\cite[Corollary 1.2]{BJ}
 \label{BJnon}
Suppose that $\Gamma<\Isom(\H^{3})$ is a finitely generated, geometrically infinite Kleinian group of the second kind such that the injectivity radius of $M=\H^{3}/ \Gamma$ is bounded away from zero. Then 
$$\dim \Lambda(\Gamma)=\dim \Lambda_{nc}(\Gamma)=2.$$
\end{theorem}

Thus, in both theorems  the Hausdorff dimension of the nonconical limit set equals the Hausdorff dimension of $\Lambda(\Gamma)$. 
In this paper we will prove that this equality does not hold in general (Theorem \ref{thm:zero}).  


The main theorems of this paper are:

\begin{theorem}
\label{thm:zero}
There is an infinitely generated discrete subgroup $\Gamma<\Isom(\H^{2})$ such that 
$$\dim \Lambda_{nc}(\Gamma)=0.$$
\end{theorem}

\begin{theorem}
\label{dim2}
Suppose that $\Gamma<\Isom(\H^{3})$ is a finitely generated, non-free, torsion-free geometrically infinite Kleinian group  such that the injectivity radius of $\H^{3}/ \Gamma$ is bounded away from $0$. Then the Hausdorff dimension of the nonconical limit set $\La_{nc}(\Ga)$ is positive. 
\end{theorem}

\begin{rem}
It is very likely that the conclusion of this theorem can be strengthened to $\dim \La_{nc}(\Ga)=2$, but proving this would require considerably more work. 
\end{rem}

The outline of the proof of Theorem \ref{dim2} is as follows: 
By Theorem \ref{BJnon}, it suffices to consider  Kleinian groups of the first kind. Moreover, using  Theorem
\ref{BJnon}, we reduce the problem to the case when $M$ is homeomorphic to the interior of a compression body which is not a handlebody; in particular, $M$ has at least two ends. 
Then there exists a finitely generated, geometrically infinite Kleinian group $\Gamma'< \Isom(\H^3)$ of the second kind such that the injectivity radius of $M'=\H^{3}/ \Gamma'$ is bounded away from $0$ and a geometrically infinite end of $M'$ is bi-Lipschitz to 
the given (geometrically infinite) end $e$ of $M$. This bi-Lipschitz homeomorphism induces a 
bi-H\"older homeomorphism from the nonconical limit set of $\Gamma'$ to set of {\em end-limit points} $\La(e)\subset \La_{nc}(\Ga)$ of the end $e$;  consequently, $\La_{nc}(\Ga)$ has positive Hausdorff dimension. 
For details, see Section \ref{maintheo}. \\

{\bf Acknowledgements.} We would like to thank Christopher Bishop and Subhadip Dey for the helpful discussions and suggestions.   During the work on this paper the first author was partly supported by the NSF grant  DMS-16-04241. The second author is grateful to Max Planck Institute for Mathematics in Bonn for its hospitality and financial support.

\section{Background}\label{sec:basics} 

\subsection{Metric geometry}

We will adopt the Bourbaki convention that neighborhoods in topological spaces need not be open. 
Let $(Y, d)$ be a metric space.  
For a subset $A\subset Y$ and a point $y\in Y$, we will denote by $d(y, A)$ the \emph{minimal distance} from $y$ to $A$, i.e.
$$d(y, A):=\inf \{ d(y, a) \mid a\in A\}.$$
Similarly, for two subsets $A, B\subset Y$ define their minimal distance as 
$$d(A, B):= \inf\{ d(a, b) \mid a\in A, b\in B \}.$$
We let $B(c,r)$ denote the open ball of radius $r$ and center $c$ in $(Y, d)$. 
We use the notation $\bar{N}_{r}(A)$ and $N_r(A)$ respectively for the closed and open $r$-neighborhoods of $A$ in $Y$: 
$$
\bar{N}_{r}(A)=\{ y\in Y: d(y, A)\leq r\}, \quad  N_{r}(A)=\{ y\in Y: d(y, A)< r\}.$$
Such neighborhoods are called {\em metric} neighborhoods of $A$. 

The Hausdorff distance $\textup{hd}(Q_{1}, Q_{2})$ between two closed subsets $Q_{1}, Q_{2}$ of $(Y, d)$ is the infimum of $r\in [0, \infty)$ such that $Q_{1}\subset \bar{N}_{r}(Q_{2})$ and $Q_{2}\subset \bar{N}_{r}(Q_{1})$. 

Given points $a, b$ in a geodesic metric space $Y$, we use the notation $ab$ for a geodesic segment from $a$ to $b$ in $Y$. Given points $a, b, c\in Y$ we let 
$[abc]$ denote a geodesic triangle in $Y$ which is the union of geodesic segments $ab, bc, ca$. 

\subsection{Ends of spaces} 

Let $Z$ be a locally path-connected, locally compact, Hausdorff topological space. The {\em ends} of $Z$ are defined as follows (see \cite{DK} for details). 
Consider an exhaustion $(K_i)$ of $Z$ by an increasing sequence of   compact subsets:
$$K_{i}\subset K_{j}, \quad \textup{ whenever } i\leq j,$$
and 
$$\bigcup_{i\in \mathbb{N}} K_{i}=Z.$$
Set $K_i^c:= Z\setminus K_i$. The ends of $Z$ are equivalence classes of decreasing sequences of connected components 
$(C_i)$ of ${K_i}^c$: 
$$C_{1}\supset C_{2}\supset C_{3}\supset \cdots $$
Two sequences $(C_i), (C'_j)$ of components of $({K_i}^c), ({K'_j}^c)$ are said to be equivalent if each $C_{i}$ contains some $C'_{j}$ and vice-versa. Then the equivalence class of a sequence $(C_i)$ is an {\em end} $e$ of $Z$. Each $C_i$ and its closure is called a 
 \emph{neighborhood} of $e$ in $Z$. The set of ends of $Z$ is denoted $Ends(Z)$. 
 An end $e$ is called {\em isolated} if it admits a closed 1-ended neighborhood $C$; such a neighborhood is called {\em isolating}.  Equivalently: There is a natural topology on the union $\hat{Z}=Z\cup Ends(Z)$  which is a compactification of $Z$ and the neighborhoods $C$ of ends $e$ as above are intersections of $Z$ with neighborhoods of  $e$ in $\hat{Z}$. Then an end $e$ is isolated if and only if it is an isolated point of $\hat{Z}$.  
 A closed neighborhood $C$ of $e$ in $Z$ is isolating if and only if $C\cup \{e\}$ is closed in  $\hat{Z}$. 
 
  A proper continuous map (a \emph{ray}) $\rho: \mathbb{R}_{+}\rightarrow  Z$ is said to be \emph{asymptotic to the end $e$} if for every neighborhood $C_{i}$ of $e$, the subset $\rho^{-1}(C_{i}) \subset \mathbb{R}_{+}$ is unbounded. One verifies that every ray is asymptotic to exactly one end of $Z$, see e.g. \cite{DK}.

\subsection{Hyperbolic spaces} 

Next, we review the notion of Gromov-hyperbolic spaces and  visual metrics on their ideal boundaries; we refer the reader to \cite{BH, BS, BS2, DK} for details. 

\begin{definition}
 Given three points $x, y, w$ in a metric space $X$, the \emph{Gromov product} of $x, y$ with respect to the 
 basepoint $w$ is defined as 
 $$
 (x \mid y)_{w}=\dfrac{1}{2}(d(x, w)+d(y, w)-d(x, y)).$$
 \end{definition}

The Gromov product is always nonnegative because of the triangle inequality. Given two real numbers $a, b\in \bar{\mathbb{R}}$, let $a\wedge b$ denote the minimum. The space $X$ is {\em $\delta$-hyperblic in Gromov's sense}, if  all $x, y, z, w\in X$ satisfy the inequality, 
$$(x\mid z)_{w}\geq (x\mid y)_{w}\wedge (y\mid z)_{w}-\delta.$$ 
A metric space $X$ is said to be \emph{Gromov-hyperbolic} if it 
is $\delta$-hyperbolic for some $\delta<\infty$. 

For {\em geodesic metric spaces}, there is an equivalent definition for Gromov hyperbolicity in terms of thinness of triangles in $X$:  A geodesic metric space is 
{\em $\eps$-hyperbolic} (in the Rips' sense) if every geodesic triangle $[xyz]\subset X$ is {\em $\eps$-slim}: 
$$
xy\subset \bar{N}_{\eps}(yz \cup zx). 
$$
This is the most common definition used in the literature and we will refer to such spaces simply as $\delta$-hyperbolic. 

We next review the notion of the {\em ideal boundary}\footnote{also known as the {\em visual} or {\em Gromov} boundary}  $\geo X$ for a $\delta$-hyperbolic (in Gromov's sense) space $X$. 

Fix a base-point $w\in X$. A sequence $\{x_{i}\}\subset X$ is said to \emph{converge at infinity} if 
$$\lim_{i, j\rightarrow \infty} (x_{i} \mid x_{j})_{w}=\infty.$$
This implies that $d(x_{j}, w)\rightarrow \infty$ as $j\rightarrow \infty$. Two sequences $\{x_{i}\}, \{y_{i}\} $ that converge at infinity are declared to be equivalent if 
$$\lim_{i\rightarrow \infty} (x_{i} \mid y_{i})_{w}=\infty.$$
This equivalence relation between sequences that converge at infinity does not depend on the choice of the base-point $w$. The \emph{ideal boundary} $\geo X$ of $X$ is defined to be set of equivalence classes of sequences that converge at infinity. For {\em proper} Gromov-hyperbolic (in Rips' sense) the ideal boundary can be defined as the set of equivalence classes of geodesic rays in $X$ where two rays are equivalent if they are Hausdorff-close, i.e. are within finite Hausdorff distance from each other.

For $a, b\in \geo X$, the Gromov product $(a \mid b)_{o}$ of $a$ and $b$ with respect to the base-point $o$ is defined as:
$$(a \mid b)_{o}=\sup \{ \liminf_{i\rightarrow \infty} (x_{i} \mid y_{i})_{o} \mid \{x_{i}\}\in a, \{y_{i}\}\in b\}.$$

Given $\varepsilon>0$, there is a standard construction of the metric on $\partial_{\infty} X$. For $x, y\in \geo X, o\in X$,  
$$d_{o, \varepsilon}(x, y):=\inf \{ \sum_{i=1}^{n} e^{-\varepsilon (x_{i-1} \mid x_{i})_{o}} \}$$
where the infimum runs over all finite sequences $x=x_{0}, x_{1}, x_{2}, \cdots, x_{n}=y$ in $\geo X$.


 \begin{definition}\cite{BS, BS2}
 A metric $d$ on $\geo X$ is \emph{visual} if there are constants $K, C>1$ such that 
 $$\dfrac{K^{-(a \mid b)_{o}}}{C}\leq d(a, b) \leq C\cdot K^{-(a \mid b)_{o}}$$
 for all $a, b\in \geo X, o\in  X$. 
 \end{definition}
 
 Consider for instance the (classical) hyperbolic $n$-space $\H^n$: This space is Gromov-hyperbolic and the standard metric on the boundary sphere is visual:

\begin{lemma}\label{lem:visual} 
Consider the unit ball model of the hyperbolic space $\H^n$. Then 
the angular metric on $S^{n-1}$ is visual. 
\end{lemma}
 \proof Let $o$ be the Euclidean center of the unit ball. Consider points $a, b\in S^{n-1}$ with the angular distance $\angle ab =\theta$.  
Then $e^{-(a \mid b)_{o}}=\sin(\theta/2)$, see \cite[Chapter 8]{roe}. It is easy to see that 
 $$\sin(\theta/2)\leq  \dfrac{\theta}{2}\leq \pi \sin(\theta/2)$$
 since $0\leq \theta\leq \pi$.
 Hence, the angular metric on the sphere $S^{n-1}$ is visual.  \qed

\medskip
For general $\delta$--hyperbolic spaces $X$, a source of visual metrics on $\geo X$ comes from: 
 
 \begin{lemma}\cite{BS}
 There is some constant $\epsilon_{0}>0$ with the following property. If $X$ is $\delta$-hyperbolic and $\varepsilon \delta\leq \epsilon_{0}$, then 
 $$\dfrac{1}{2} e^{-\varepsilon (a \mid b)_{o}}\leq d_{o, \varepsilon}(a, b)\leq e^{-\varepsilon (a \mid b)_{o}}, \quad \forall a, b\in \geo X.$$
 \end{lemma}

\medskip 
 Thus, for a $\delta$-hyperbolic space $X$ and every $w\in X$ there exist $\varepsilon>0$ such that $d_{o, \varepsilon}$ is a visual metric on $\geo X$. One important property we will be using to prove Theorem \ref{dim2} is that quasi-isometries between Gromov hyperbolic spaces induce bi-H\"older maps on the ideal boundaries in terms of the visual metrics. 
 
 \begin{definition}
A map $f: X\rightarrow Y$ between metric spaces $(X, d_{1}), (Y, d_{2})$ is said to be a $(\kappa,c)$-\emph{quasiisometric embedding} if 
$$\dfrac{1}{\kappa}d_{1}(x, x')-c\leq d_{2}(f(x), f(x'))\leq \kappa d_{1}(x, x')+c$$
for all $x, x'\in X$. A quasiisimetric embedding $f$ is called a {\em quasiisometry} if, in addition,  one has 
$\textup{hd}(f(X), Y)\le c$. 
 \end{definition}
 
A \emph{quasi-geodesic} in a metric space $X$ is a quasiisometric map $\gamma: I\rightarrow X$ where $I\subset \mathbb{R}$ is an interval. 

\begin{theorem}
[Hyperbolic Morse Lemma, see e.g. \cite{DK}] Let $X$ be a $\delta$-hyperbolic geodesic metric space. 
Then the image of any $(\kappa,c)$-quasigeodesic in $X$ is within distance $D(\kappa,c, \delta)$ 
from a geodesic in $X$ with the same end-points. 
\end{theorem}

\begin{theorem}[See e.g. \cite{BS}]
\label{holder}
Suppose that $X$ and $Y$ are $\delta$-hyperbolic geodesic spaces. If $f: X\rightarrow Y$ is a $(\kappa, c)$-quasiisometry (resp. quasiisometric embedding), then $f$ induces a bi-H\"older homeomorphism (resp. a  bi-H\"older embedding) 
between boundaries $\geo f: \geo X \rightarrow \geo Y$ equipped with visual metrics $d_{\geo X}, d_{\geo Y}$, more precisely: 
$$
C_{2} (d_{\geo X} (a, b))^{1/\kappa}\leq  d_{\geo Y}(\geo f(a), \geo f(b))\leq C_{1} (d_{\geo X}(a, b))^{\kappa}$$
where $C_{1}, C_{2}$ are constants depending on $\delta, \kappa$ and $c$. 
\end{theorem}

\subsection{Hausdorff dimension}
 
We now  review the Hausdorff dimension of metric spaces. 
For $s>0$, the $s$-dimensional \emph{Hausdorff content} of a metric space  $E$ 
is defined as:
$$
\mathcal{H}^{s}(E)= \lim_{\delta\to 0} \inf \{ \sum_j r^s_{j} \},$$
where the infimum is taken over all covers $\{B(x_{j}, r_{j}): j\in J\}$ of $E$ by metric balls  satisfying
$$
\sup_{j} r_j\le \delta. 
$$


\begin{definition}
The \emph{Hausdorff dimension} of $E$ is defined as
$$
\dim(E): =\inf\{\alpha\in \R_+: \mathcal{H}^{\alpha}(E)=0\}=\sup \{ \alpha:  \mathcal{H}^{\alpha}(E)=\infty\}.$$
\end{definition}

In particular, if $\mathcal{H}^{\alpha}(E)<\infty$ then $\dim(E)\le \al$.

\subsection{Kleinian groups}

We  turn to  the classical hyperbolic space $\H^{n}$ with the visual ideal boundary $S^{n-1}=\geo \H^n$. It is a uniquely geodesic space; we will use the notation $\delta_0$ for a hyperbolicity constant on $\H^n$, i.e. for the constant such that every geodesic triangle in $\H^n$ is $\delta_0$-slim. 
(One can take $\delta_0=2 \cosh^{-1}(\sqrt{2})$.) We let $xy\subset \H^{n}$ denote the geodesic segment connecting $x, y\in \H^{n}$. Similarly, given $x\in \H^{n}$ and $\xi\in S^{n-1}$ we use the notation $x\xi$ for the unique geodesic ray emanating from $x$ asymptotic to $\xi$; for two distinct points $\xi, \eta\in S^{n-1}$, we use the notation $\xi\eta$ to denote the  geodesic asymptotic to $\xi$ and $\eta$; it is unique up to reparameterization.  The {\em closed convex hull} $C(\Lambda)$ of a subset $\La\subset S^{n-1}$ is the smallest closed convex subset of $\H^{n}$ whose accumulation set in $\H^{n}\cup S^{n-1}$ equals $\Lambda$. This set exists whenever $\La$ has cardinality $\ge 2$. 

For a pair of points $x, y\in \H^{n}$, we let $H(x, y)$ denote the closed half space in $\H^{n}$ given by 
$$H(x, y)=\{p\in \H^{n}: d(p, x)\leq d(p, y) \}.$$ 
Given a Kleinian group $\Gamma<\Isom(\H^{n})$ and a point $p\in \H^n$ with trivial $\Ga$-stabilizer, 
the \emph{Dirichlet fundamental domain for $\Gamma$ centered at $p$ } is defined to be the set 
$$D_{p}(\Gamma)=\bigcap_{\gamma\in \Gamma-\{1\}} H(p, \gamma(p)). $$

Recall that the limit set $\La=\Lambda(\Gamma)$ of a discrete subgroup $\Ga< \Isom(\H^n)$ has cardinality $0, 1, 2$ or continuum; 
a subgroup $\Gamma$ is called \emph{elementary} if $\Lambda(\Gamma)$ is finite. Otherwise, it is \emph{nonelementary}. In this paper, we concentrate on nonelementary Kleinian groups $\Gamma$. For a nonelementary Kleinian group $\Ga$, $C(\La)/\Ga$ is a convex subset of $M=\H^n/\Ga$, called the {\em convex core} $C(M)$ of  $M$. A nonelementary Kleinian group is called {\em convex-cocompact} if $C(\La)/\Ga$ is compact.


\medskip 


A limit point $\la\in \La(\Ga)$ of a Kleinian group $\Ga< \Isom(\H^n)$ 
is called {\em nonconical} if the projection of one (equivalently, every) geodesic ray $x\la\subset \H^n$ to 
$M=\H^n/\Ga$ is a proper map. The set of nonconical limit points of $\Ga$ is denoted $\Lambda_{nc}$. Its complement $\La(\Ga) \setminus \La_{nc}$ is the {\em conical limit set} $\La_c$ of $\Ga$.  A Kleinian subgroup $\Ga$ is convex-cocompact if and only if $\La(\Ga)=\La_c(\Ga)$.

\medskip 
{\bf Critical exponent.} The  \emph{critical exponent} (or the \emph{Poinc\'are exponent} or the {\em exponent of convergence}) 
of a Kleinian group $\Gamma< \Isom(\H^n)$ is defined as:  
$$
\delta(\Gamma):=\inf\{ s: \sum_{\gamma\in \Gamma} \exp(-sd(p, \gamma(p)))<\infty \},$$
where $p\in \H^{n}$. The critical exponent depends only on $\Ga$ and not on $p$.  
The critical exponent of a nonelementary Kleinian group equals the Hausdorff dimension of its conical limit set:
$$
\delta(\Ga)= \dim ( \La_{c}(\Ga)),
$$
see \cite{BJ}.

\medskip 
{\bf Klein combination.} Combination theorems provide a useful procedure for constructing  Kleinian  groups. Suppose that $\Ga_1, \Ga_2< \Isom(\H^n)$ are Kleinian subgroups, and $\geo \H^n$ is expressed as a union of compact subsets, $F_1\cup F_2$, such that:
$$
\gamma F_i\cap F_i=\emptyset, ~~ \forall \gamma\in \Ga_i \setminus \{1\}, i=1, 2. 
$$
Then the pair of subgroups $\Ga_1, \Ga_2$ is said to satisfy the conditions of the {\em Klein combination theorem}. Under these conditions one has:

\begin{theorem}\label{thm:combination} 
The subgroup $\Ga<  \Isom(\H^n)$ generated by $\Ga_1, \Ga_2$ is again Kleinian and is naturally isomorphic to the free product $\Ga_1\star \Ga_2$. Furthermore, $\La_{nc}(\Ga)$ is the $\Ga$-orbit of $\La_{nc}(\Ga_1)\cup \La_{nc}(\Ga_2)$. 
\end{theorem}

See for instance \cite[Theorem C.2, section VII.C]{Mas}. Note that Maskit states and proves this theorem only for $n=3$ (we will need it only for $n=2$), but the proof is general.  

Not every pair of Kleinian subgroups satisfies the conditions of the Klein combination theorem: For instance, these conditions imply that $F_i\subset \Omega(\Ga_i)$, $i=1,2$ and, hence, $\Lambda(\Ga_1)\cap \Lambda(\Ga_2)=\emptyset$. However, if $\Gamma_i< \Isom(\H^n)$ are Kleinian with $\Omega(\Gamma_i)\ne \emptyset, i=1, 2$, then there exist conjugates $\Gamma'_i, i=1, 2$, of $\Ga_1, \Gamma_2$, which do satisfy these conditions. 

\medskip
\subsection{Ends of hyperbolic 3-manifolds}\label{sec:ends-3D} 

Suppose that $M$ is a complete connected hyperbolic 3-manifold with finitely generated fundamental group. According to the solution of the ELC (Ending Lamination. Conjecture) in the work of Minsky \cite{Minsky} and Brock--Canary--Minsky \cite{BCM} (see also alternative proofs by Soma \cite{Soma2} and Bowditch \cite{Bowditch:ELC}), the geometry of $M$ is completely determined by its topology and a certain set of {\em asymptotic invariants} of the {\em ends } of $M$. For simplicity, we only discuss this in the case when $M$ has no cusps. 
According to \cite{Agol, Bonahon, CG} (see also Soma's paper \cite{Soma1}), 
the manifold $M$ is {\em topologically tame}, i.e. is homeomorphic to the interior of a compact manifold with boundary $\bar{M}$.  
Thus, ends $e_1,...,e_k$ of $M$ are in bijective correspondence with the boundary surfaces $S_1,...,S_k$ of $\bar{M}$. Each end $e_i$ of $M$ is either {\em geometrically finite} or {\em geometrically infinite}. An end $e$ of $M$ is \emph{geometrically finite} if it has an isolating neighborhood $E$  disjoint  from 
the convex core $C(M)$; otherwise, it is \emph{geometrically infinite}. Each geometrically infinite end $e$ is {\em simply degenerate}, i.e. 
has a closed isolating neighborhood $E$ homeomorphic to $S\times [0, 1)$ (where $S$ is a compact surface) and there exists a sequence of {\em pleated surfaces} $S_{n}$ in $E$ leaving every compact set such that for each $n$, $S_{n}$ is homotopic to $S\times \{0\}$ within $E$. We refer the reader to 
 \cite{Bonahon, C1} for more detail.

Let $M_c$ denote the convex core of $M$. Let $E_1,...,E_l$ denote closures of the connected components of $M -M_c$; topologically speaking, these are products $S_i\times (0,1]$, $i=1,...,l$. The subsets $E_i, i=1,...,l$, serve as closed isolating neighborhoods of geometrically finite ends of $M$. 
Thus, we index the ends of $M$ so that the ends $e_1,...,e_l$ are geometrically finite. The end-invariants of the geometrically finite ends are (marked) Riemann surfaces $X_1,...,X_l$ defined as follows. Let $\tilde{E}_i$ denote a component of the preimage of $E_i$ in $\H^3$ (with respect to the universal covering map $\H^3\to M$). Each proper geodesic ray $\rho: \R_+\to E_i$ lifts (nonuniquely) to a geodesic ray $\tilde\rho$ in  $\tilde{E}_i$. 
The set of limit points $\tilde\rho(\infty)$ of these lifts forms an open subset $\Omega_i$ of $S^2=\geo \H^3$. The group $\Ga_i$ of covering transformations of $\tilde E_i\to E_i$ acts properly discontinuously on $\Omega_i$ and the Riemann surface $X_i$ equals $\Omega_i/\Gamma_i$. The neighborhoods $E_i$ of the end $e_i$ then admits a compactification 
$$
\bar{E}_i= (\tilde{E}_i \cup \Om_i)/\Ga_i.
$$
An alternative description of the Riemann surface $X_i$ is given as the limit (as $\kappa\to 0$) of conformal structures of convex surfaces $S_i(\kappa)$ of constant  curvature $\kappa\in (0,-1)$ foliating $E_i= S_i\times (0,1]$, see \cite{Labourie}. 

While asymptotic invariants parameterizing geometrically finite ends of $M$ belong to the Teichm\"uller spaces $T(S_i)$ of surfaces $S_i$, $i=1,...,l$,  asymptotic invariants parameterizing the geometrically infinite ends belong to quotient spaces of Thurston boundaries of $T(S_i)$'s $i=l+1,...,k$. 
More precisely, for each geometrically infinite end $e_i$, $i=l+1,...,k$,  there exists a sequence $\alpha_n$ of simple essential loops on $S_i$ such that the corresponding sequence of closed geodesics $\al^*_n$ in $M$ is contained in $E_i$ and {\em escapes} the end $e_i$, i.e. every compact subset $K\subset E_i$ intersects only finitely many members of the sequence $(\al^*_n)$. The sequence of loops $(\al_n)$ defines a sequence 
$(\bar\al_n)$ in the {\em space of projective classes of measured geodesic laminations on $S_i$}, $PML(S_i)$. Here we equip $S_i$ with some background hyperbolic metric. (The union $T(S_i)\cup PML(S_i)$ admits a certain natural topology making it a closed ball compactification of $T(S_i)$.) The space $PML(S_i)$ is compact and one considers the set $L_i\subset PML(S_i)$ of accumulation points of sequences   $(\bar\al_n)$. (In fact, it suffices to take just one sequence.) Any two elements of $L_i$ are represented by measured geodesic laminations on $S_i$ which differ  only by the transverse measure and, thus, the transverse measures have the same support sets $\epsilon_i$. The geodesic lamination $\epsilon_i$ is called the {\em ending lamination} of the end $e_i$. 

Thus, one obtains the set of asymptotic end-invariants $(X_1,...,X_l, \eps_{l+1},...,\eps_k)$ of the manifold $M$. According to the ELC, the manifold $M$ is uniquely determined by its topology and its set of asymptotic invariants. Furthermore, this uniqueness theorem has an existence counterpart: Given a compact topological 3-manifold $\bar{M}$, under certain conditions on the topology of $M$ and on the end-invariants associated to the boundary surfaces of $\bar{M}$, 
they can be realized as end-invariants of a certain complete hyperbolic structure on the interior of $M$,  
see the paper by Namazi, Souto and Ohshika \cite{NS, Ohshika1, Ohshika2}. We will state only a weak form of this existence result: 

\begin{thm}
Given a hyperbolic manifold $M$ all whose ends are geometrically infinite with the end-invariants $(\eps_{1},...,\eps_k)$ and given any collection $(Y_1,...,Y_l), 1\le l\le k$, of marked Riemann surface structures on the boundary surfaces $S_1,...,S_l$   of $\bar{M}$, there exists a hyperbolic structure on $M$ with the end-invariants 
$$
(Y_1,...,Y_l, \eps_{l+1},..., \eps_k).  
$$
\end{thm}

\subsection{Types of nonconical limit points} 

Suppose that $e$ is an end of a complete hyperbolic manifold $M=\H^n/\Ga$. The {\em limit set} $\La(e)$ of $e$ is the subset of $\La_{nc}(\Ga)\subset \La=\La(\Ga)$ consisting of 
limit points $\la\in \La(\Ga)$ such  that for some (equivalently, every) $x\in \H^n$ the geodesic ray $x\la$ in $\H^n$ projects to a proper ray $\rho$ in $M$ asymptotic to the end $e$. While $\rho$ is proper, it can diverge to infinity in $M$ in different ways. The limit points $\la\in \La(e)$ are classified accordingly as: 

\begin{definition}
For $\beta\in (0, 1]$ a limit point $\la$ is called $\beta$-{\em deep} for  if 
$$
\lim \inf_{t\to \infty} \frac{d(\rho(0), \rho(t)) }{t} \ge \beta. 
$$
\end{definition}


A {\em J{\o}rgensen limit point} is a  point $\xi\in \Lambda(\Gamma)$ such that there exists a geodesic ray $x\xi$ 
asymptotic to $\xi$ which is completely contained in some Dirichlet domain $D_{p}(\Gamma)$ of $\Ga$. 
This definition is easily seen to be equivalent to the condition that the projection of 
the ray $p\xi$ to $M$ is an isometric embedding, i.e. is a geodesic in the sense of metric geometry. Thus, J{\o}rgensen limit points are $1$-deep. Let $\Lambda_{\beta}(\Gamma)$ denote the set of limit points which are $\beta$-deep. Bishop \cite{Bi2} proved that 
$$\dim \La_{nc}(\Gamma)=\dim(\bigcup_{0<\beta \leq 1} \Lambda_{\beta}(\Gamma)).$$
This result was sharpened by Gonye in \cite{Gonye}. 

In the paper we will need a variation on the notion of deep limit points which is neither weaker nor stronger than the one given above. We assume that $e$ is an isolated geometrically infinite end of a complete hyperbolic manifold $M=\H^n/\Ga$ and let $E\subset M$ be its isolating neighborhood.

\begin{lemma}\label{lem:special ray} 
Every closed isolating neighborhood $E$ of $e$ contains an isometrically embedded geodesic ray $\rho: \R_+\to M$ 
(necessarily asymptotic to $e$) such that for every $t\ge 0$ 
$$
t= d(\rho(t), \partial E). 
$$
\end{lemma}
\proof 
For each $i\in \N$ we let $z_i w_i$ denote a shortest geodesic (necessarily of length $i$) in $M$ from $\partial E$ to the $i$-level set  of the distance function $d(\cdot, \partial E)$ on ${E}$. Since the sequence $(z_i)$ 
lies in the compact $\partial E$, the sequence of geodesics  $z_i w_i$ subconverges to a geodesic ray $\rho$ in $M$ contained entirely in ${E}$. 
By the construction, the limiting ray $\rho$ satisfies the desired properties. \qed  

\begin{lemma}
Assuming, in addition, that $n=3$ and the end $e$ of $M=\H^{3}/\Gamma $ is simply degenerate, then the geodesic ray $\rho$ in Lemma \ref{lem:special ray} can be taken to be equal to the projection of some ray $x\lambda$ where $x\in \H^{3}$ and $\la\in \La(e)$. 
\end{lemma}

\proof
By the definition of simply degenerate end, we assume that the neighborhood $E$ is contained in  the convex core $C(M)$. Then the preimage of $E$ in $\H^{3}$ is contained in the convex hull $C(\Gamma)$, and there is a lift $\tilde{\rho}=x\la$ of $\rho$ also contained in $C(\La), \La=\La(\Ga)$. 
Note that  $\partial_{\infty} C(\La)=\Lambda$. Hence,  $\la\in \Lambda$, i.e. $\la\in \Lambda(e)$. 
\qed

\medskip 

We fix a connected component $\tilde E$ of the preimage of $E$ in $\H^n$. Then the end-limit set $\La(e)$ is contained in the union of $\Ga$-translates of the 
accumulation set $\geo \tilde E$ in $S^{n-1}= \geo \H^n$. Fix a point $x\in \tilde E$ and a constant $D>0$. 

\begin{definition}
A limit point $\la\in \La(e)\cap \geo \tilde E$ is {\em $(x,D)$-deep} (with respect to $E$) if the  geodesic ray $x\la$ is disjoint from the $D$-neighborhood of 
$\partial \tilde E$. 
We let $\La_{(x,D)}(E)$ denote the subset of $\La(e)$ consisting of $(x,D)$-deep limit points of the end $e$. 
\end{definition}

Let $\Ga_E$ denote the stabilizer of $\tilde E$ in $\Ga$. 
Thus, for $\ga\in \Ga_E$, $\la$ is $(x,D)$-deep if and only if $\ga(\la)$ is $(\ga(x), D)$-deep. 

We next 
consider   
$E$ as a metric space with the distance function $d_E$ obtained by restricting the Riemannian 
distance function of $M$. Since ${E}$ is one-ended, the quasiisometry class of the metric space ${E}$ is independent of the choice of $E$. 
While in general, the coarse geometry of ${E}$ can be quite complicated, we will be primarily interested in the case when $E$ is quasiisometric to the half-line. This will be the case for every end of a tame hyperbolic 3-manifold with finitely-generated fundamental group and injectivity radius bounded below, see 
e.g. \cite{C1}. 

Given $E$, let $E_D$ denote the  $D$-level set  of the distance function $d(\cdot, \partial E)$ on ${E}$.

\begin{lemma} 
The following conditions are equivalent for a closed isolating neighborhood $E$ of an isolated end $e$ of a hyperbolic manifold $M$. 

1. $({E}, d_E)$ is quasiisometric to $\R_+$. 

2. For every geodesic ray $\rho$ in $M$ asymptotic to $E$, a metric neighborhood of $\rho$ contains $E$. 

3. For an isometrically embedded geodesic ray $\rho$ in $E$ as in Lemma \ref{lem:special ray}, a metric neighborhood of $\rho$ contains $E$.

4. There exists a constant $a$ such that for every $D$, the diameter (with respect to $d_E$) of $E_D$  is $\le a$. 
\end{lemma}
\proof (1) implies (2)  since every continuous proper map $\R_+\to \R_+$ is coarsely surjective. The implications (2) $\Rightarrow$ (3) and (3) 
 $\Rightarrow$ (1) are immediate. (For the last implication, we note that the (coarsely well-defined) nearest-point projection $P: E\to \rho(\R_+)$ is the coarse inverse of the isometric embedding $\R_+\to E$ given by $\rho$.)

(3) $\Rightarrow$ (4). Let $\rho$ be an isometrically embedded geodesic ray in $E$  as in Lemma \ref{lem:special ray} and let 
$P: E\to \rho(\R_+)$ be a  the nearest-point projection. Set $R_0:= \diam(\partial E)$ and let $R_1$ be an upper bound on the diameters of point-preimages under 
$P$. 

For $z\in E_D$, its projection $P(z)=\rho(t)$ satisfies the inequalities
$$
d(P(z), z)\le R_1, \quad D- (R_0+R_1)\le t\le D + (R_0+R_1).
$$
In particular, $d(P(z), \rho(D))\le R_0+R_1$. Since $d(\rho(D),\partial E)=D$, 
 the diameter of $E_D$ is at most $2(R_0+2R_1)$. 

(4) $\Rightarrow$ (3). Every $z\in E_D$ the distance from $z$ to $\rho(D)$ is $\le a$. Hence, $\bar{E}$ is contained in the closed $a$-neighborhood 
of the image of $\rho$. \qed 

\begin{definition}
When one of the conditions in this lemma holds, we will say that the end $e$ is {\em narrow}. 
In the case when (4) holds for a specific constant $a$, we will say that  $E$ is {\em $a$-narrow}. 
\end{definition} 

\begin{lemma}\label{lem:limit-of-narrow-end} 
If $e$ is a narrow end then for every geodesic ray $\rho$ in $M$ asymptotic to $E$, each lift $\tilde\rho=x\xi$ of $\rho$ to $\H^n$ yields a limit point $\xi$ of $\Ga_E$.  
\end{lemma}
\proof If $\xi$ is not a limit point of $\Ga_E$ then it is a point in the domain of discontinuity of this group, which implies that the injectivity radius of $M$ at $\rho(t)$ diverges to infinity as $t\to \infty$. This contradicts the narrowness assumption. \qed

\begin{lemma}
Suppose that $E$ is $a$-narrow and $x\in \tilde E$ is such that $d(x, \partial \tilde E)>  a+ \delta_0$. Then for $D:= d(x, \partial \tilde E)- (a+\delta_0)$, we have 
$$
\La(e)\cap \geo \tilde E= \bigcup_{\ga\in \Ga_E} \ga \La_{(x, D)}(E). 
$$
\end{lemma}
\proof Since one inclusion is clear, we need to prove that 
$$
\La(e)\cap \geo \tilde E \subset \bigcup_{\ga\in \Ga_E} \ga \La_{(x, D)}(E). 
$$
Let $z\in M$ denote the projection of $x$. Given $\la\in \La(e)\cap \geo \tilde E$ and the ray $x\la$, we let $\rho$ denote the projection 
of $x\la$ to $M$: The image is a geodesic ray in $M$.  Since $\rho$ is a proper ray asymptotic to the end $e$, 
$\rho$ contains a maximal subray $\rho_1$ emanating from a point $z_1\in E_{D_{x}}$, such that 
$\rho_1\cap N_{D_{x}}(\partial E)= \{z_1\}$ where $D_{x}=d(x, \partial \tilde{E})$.
In other words, a subray $x_1 \la\subset x\la$ is $(x_1, D_{x})$-deep, where $x_1\in \tilde E$ projects to $z_1$.   Note that if $\rho\cap N_{D_{x}}=z$, then $z_{1}=z, \rho_{1}=\rho$ and $x\la$ is $(x, D)$-deep.

By the $a$-narrowness of $E$, for each pair of points $z, z_1$ as above, there exists a geodesic $s$ in $M$ of length $\le a$ 
connecting $z$ to $z_1$. The concatenation $s\star \rho_1$ is properly homotopic to a unique geodesic proper ray $r$ emanating from $z$. Of course, the rays $\rho, r$  need not be properly homotopic. The segment $s$ and the ray $r$ lift to  a geodesic segment $\tilde s$ and a geodesic ray $\tilde r$ in $\tilde E$, so that (for some $\ga\in \Ga_E$) the ray $\tilde r$ connects a point $\ga(x)$ to $\la$ and $\tilde s= \ga(x) x_1$. Since $\H^n$ is $\delta_0$-hyperbolic, $\tilde r$ is contained in the $\delta_0$-neighborhood of $\tilde \rho_1 \cup \tilde s$, i.e. in the $(\delta_0+a)$-neighborhood of $\tilde \rho_1$. Since $\tilde \rho_1$ intersects the $D_{x}$-neighborhood of $\partial \tilde E$ only at $x_{1}$ and $D=d(x, \partial \tilde E) - (a+ \delta_0)$, we conclude that $\la$  is $(\ga(x), D)$-deep. \qed 

\medskip 
Since we are only interested in Hausdorff dimensions, we conclude:

\begin{cor}\label{cor:edim} 
$$
\dim \La(e) = \dim \La_{(x, D)}(E). 
$$
\end{cor}

\section{Proofs of the main results}
\label{maintheo}


\textbf{Proof of Theorem \ref{thm:zero}}:
  We will be working with the upper half-plane model of the hyperbolic plane $\H^2$. 
Consider an infinite collection of pairwise disjoint closed Euclidean disks $B_i=B(c_i, r_i)$ centered at points $c_i\in \Z\subset \R$, indexed by positive integers $i\in \mathbb{N}$.  The semicircles $C_i$ which are the intersections of $\partial B_i$ with $\H^2$, are hyperbolic geodesics.  

We let $h_i\in \Isom(\H^2)$ be the inversion in the boundary circle of the ball $B_i$. 
Let $\Ga$ denote the subgroup of  $\Isom(\H^2)$ generated by $S=\{h_i: i\in \N\}$. 
For $k\in \N_0$ we will denote by $\Ga_k$ 
the subgroup of $\Ga$ generated by the subset $S_k=\{h_i: i > k \}$. Thus, $\Ga=\Ga_0$. 

Define 
$$
\Phi_k= \H^2 \cup \R \cap \bigcap_{i> k} \mbox{Ext} (B_i). 
$$
By the Poincar\'e Fundamental Domain Theorem, $\Ga$ is isomorphic to the free product of countably many $\Z_2$'s generated by the involutions $h_i, i\in \N$. By the same theorem, $\Ga_k$ is a discrete subgroup of $\Isom(\H^2)$ and 
$\Phi_k$ is a fundamental domain for the action of $\Ga_k$ on $\H^2 \cup \Omega(\Ga_k)$, where $\Omega(\Ga_k)$ is the discontinuity domain for the action of $\Ga_k$ on the ideal boundary of $\H^2$.

In particular, we can describe the limit set $\La_k$ of $\Ga_k$ as follows. For $n\in \N$ let $W_{k,n}$ denote the subset of $\Ga_k$ consisting of words of length $\le n$ with respect to the generating set $S_k$. Let $H_{k,n}$ denote the complement in $\H^2 \cup \R$ 
to the union of images of $\Phi_k$ under the elements of $W_{k,n}$. Thus, $H_{k,n}$ is a disjoint union of round disks and 
$$
\La_k= \{\infty\} \cup \bigcap_{n\in \N} H_{k,n}. 
$$
Since we are interested only in the Hausdorff dimension, 
\begin{equation}\label{eq:intersection} 
\dim \La_k= \dim  \bigcap_{n\in \N} H_{k,n}. 
\end{equation}

We also observe that for each $k$, the group $\Ga$ is obtained by Klein combination of the geometrically finite Kleinian 
group $\Ga'_k$ generated by $h_1,...,h_{k}$, and the subgroup $\Ga_k$. In particular, by Theorem \ref{thm:combination}, 
$$
\dim \La_{nc}(\Ga)= \dim \La_{nc}(\Ga_k)\le \dim \La_k.  
$$

We will prove that for a suitable choice of the centers $\{c_i: i\in \mathbb{N}\}$ and the radii $\{r_{i}: i\in \mathbb{N} \}$, 
$$
\lim_{k\to\infty} \dim \La_k=0,
$$
which will imply that the Hausdorff dimension of the nonconical limit set $\La_{nc}(\Ga)$ is zero. 

We say that a multi-index $\underline{i}= (i_1,...,i_n)\in \mathbb{N}^n$ is {\em reduced} if 
$i_m \neq  i_{m+1}$ for all $m$. We let $I_{k,n}$ denote the set of reduced multi-indices which belong to $(k,\infty)^n$.  
The number $n$ is the {\em length} of $\underline{i}$, $n=\ell(\underline{i})$. For 
 every reduced multi-index we set 
$$
h_{\underline{i}}:= h_{i_1}\circ ... \circ h_{i_n}. 
$$
This is a reduced word in $S$.  We also set $C_{\underline{i}}:= h_{i_{1}}\circ \cdots \circ h_{i_{n-1}} (C_{i_{n}})$. We let $c_{\underline{i}}$ and $r_{\underline{i}}$ denote the center and the radius of the semicircle $C_{\underline{i}}$.  For $\underline{i}\in I_{k,n}$, $C_{\underline{i}}$ is a boundary geodesic of $H_{k,n}$. 

Since $\Phi$ is a fundamental domain of $\Ga$, it follows that for each compact subset $K\subset \R$ and a sequence 
of multi-indices  $\underline{i}_m$ satisfying 
$$
c_{\underline{i}_m}\in K, \quad \lim_{m\to\infty} \ell(\underline{i}_m)=\infty,
$$
we have
\begin{equation}\label{eq:r-to-0}
 \lim_{m\to\infty} r_{\underline{i}_m}=0.
\end{equation}

Define 
$$
\mu:= \sup_{i\ne j} \frac{1}{|c_i -c_j| -1}. 
$$

\begin{lemma} [A. Beardon; Lemma 1 in \cite{Be}] 
\label{induction}
If $r_{i}\leq 1$ for every $i\in \mathbb{N}$ then every reduced multi-index $\underline{i}= (i_1,...,i_n)$ satisfies the inequality 
\begin{equation}
\label{radius induction}
r_{i_{1} \cdots i_{n}} \leq   \dfrac{r_{i_{2} \cdots i_{n}}}{(|c_{i_{1}}-c_{i_{2}}|-1)^{2}}\leq    \mu^{2}r_{i_{2} \cdots i_{n}}
\end{equation}
\end{lemma}


\begin{lemma}
\label{dim-smaller3}
Suppose that 
\begin{equation}\label{eq:al-sum} 
\sum_{i=k+1}^\infty r_{i}^{\alpha}<\infty
\end{equation} 
 and for each $n$, 
\begin{equation}
\label{dim-smaller}
\sum_{\underline{i}=(i_{1}, \cdots, i_{n})} r_{\underline{i}}^{\alpha} \leq  \sum_{\underline{j}=(j_{1}, \cdots, j_{n-1})} r_{\underline{j}}^{\alpha},
\end{equation}
where the sums are taken over $I_{k,n}$ and $I_{k,n-1}$ respectively. 
Then $\dim(\Lambda_k)\leq \alpha$. 
\end{lemma}
\proof It suffices to show that for every compact $K\subset \La_k\cap \R$,
\begin{equation}\label{eq:<infty} 
{\mathcal H}^\al(K)<\infty. 
\end{equation} 
By \eqref{dim-smaller} and \eqref{eq:al-sum}, for every $n$ we have 
$$
\sum_{\underline{i}=(i_{1}, \cdots, i_{n})} r_{\underline{i}}^{\alpha} \leq \sum_{i=k+1}^\infty r_{i}^{\alpha}<\infty. 
$$
We restrict this sum to those multi-indices $\underline{i}$ 
for which $B(x_{\underline{i}}, r_{\underline{i}})\cap K\ne \emptyset$. 
In view of \eqref{eq:r-to-0}, for these multi-indices $\underline{i}$ 
$$
\ell(\underline{i})\to\infty \Rightarrow  r_{\underline{i}}\to 0. 
$$
Thus, \eqref{eq:<infty} holds. \qed

\begin{lemma}
The inequality 
\begin{equation}
 \label{centercontrol}
\sum_{i\neq j, i> k, j> k} \left( \dfrac{1}{|c_{i}-c_{j}|-1}\right)^{2\alpha}\leq 1
\end{equation}
implies \eqref{dim-smaller}. 
\end{lemma} 
\proof  Observe that 
\begin{equation}
\begin{aligned}
\sum_{\underline{i} =(i_{1}, \cdots, i_{n})\in I_{k, n}} r^{\alpha}_{\underline{i}} &=\sum_{ i_{1}\neq i_{2}, i_{1}>k} \left(\sum_{(i_{2}, \cdots, i_{n})\in I_{k, n-1}} r_{i_{1} \cdots i_{n}}^{\alpha}\right)\\
& \leq \sum_{i_{1}\neq i_{2}, i_{1}>k}  \left( \sum_{(i_{2}, \cdots, i_{n})\in I_{k, n-1}} \dfrac{r_{i_{2}\cdots i_{n}}^{\alpha}}{(|c_{i_{1}}-c_{i_{2}}|-1)^{2\alpha}}\right),
\end{aligned}
\end{equation}
where the inequality comes from \eqref{radius induction}. By \eqref{centercontrol},
$$\sum_{i_{1}\neq i_{2}, i_{1}>k} \dfrac{1}{(|c_{i_{1}}-c_{i_{2}}|-1)^{2\alpha}}\leq 1 $$
for each fixed $i_{2}>k$. Hence,
$$\sum_{i_{1}\neq i_{2}, i_{1}>k} \left( \sum_{(i_{2}, \cdots, i_{n})\in I_{k, n-1}} \dfrac{r_{i_{2}\cdots i_{n}}^{\alpha}}{(|c_{i_{1}}-c_{i_{2}}|-1)^{2\alpha}}\right)\leq \sum_{(i_{2}, \cdots, i_{n})\in I_{k, n-1}} r^{\alpha}_{i_{2}\cdots i_{n}},$$
which implies \eqref{dim-smaller}. \qed

\bigskip 
We now choose $c_i$'s and $r_i$'s so that \eqref{eq:al-sum} and \eqref{centercontrol} hold  for every 
$k\ge 2$ with $\al=\frac{1}{2k}$.  
First, take $r_{i}:= 2^{-2i^{2}}$. Then, clearly, for every $k\geq 2$,
$$\sum_{i=k+1}^{\infty} r_{i}^{1/2k}\leq \sum_{i\in \mathbb{N}} \dfrac{1}{2^{k+2i}}=\dfrac{1}{3\cdot 2^{k}}\leq 1,$$
which implies \eqref{eq:al-sum}.

Furthermore, we set $c_{1}=0$ and inductively define 
$$
c_{i}= c_{i-1}+2^{i^{2}+2}+1.
$$ 
We now verify \eqref{centercontrol}, i.e. that  
\begin{equation}
 \label{centercontrol2}
\sum_{i\neq j, i> k, j> k} \left( \dfrac{1}{|c_{i}-c_{j}|-1}\right)^{1/k}\leq 1
\end{equation}

For a fixed $i>k$, we have  
$$
\sum_{j>i} \left(\dfrac{1}{|c_{i}-c_{j}|-1}\right)^{1/k}\leq \dfrac{1}{2^{i+2}}+\dfrac{1}{2^{i+4}}+\dfrac{1}{2^{i+6}}+\cdots\leq \dfrac{1}{3\cdot 2^{i}}.$$
Thus,  for each $k\ge 1$, 
$$
\sum_{ i\neq j, i>k, j >k}  \left(\dfrac{1}{|c_{i}-c_{j}|-1}\right)^{1/k}\leq \dfrac{2}{3\cdot 2^{k+1}}+\dfrac{2}{3\cdot 2^{k+2}}+\cdots\leq \dfrac{1}{3\cdot 2^{k-1}}\leq 1.$$
Therefore, for every $k$, $\dim \La_k \le \al= \frac{1}{2k}$, which implies that $\dim \La_{nc}(\Ga)=0$. \qed

\medskip 
\textbf{Proof of Theorem \ref{dim2}}: 
Since the group $\Gamma\cong \pi_1(M)$ is finitely generated,  the quotient manifold $M=\H^{3}/ \Gamma$ is 
{\em tame}, i.e. is diffeomorphic to the interior of a compact manifold $\bar{M}$. The ends $e_i, i=1,...,N$, 
of $M$ are in bijective correspondence with the boundary components $S_i$ of $\bar{M}$.


 The group $\Gamma$ is geometrically finite if and only if all the ends of $M$ are geometrically finite. 
 Since $\Gamma$ is assumed to be geometrically infinite, one of the ends, say, $e_1$, is geometrically infinite. More precisely, 
 there exists a sequence of pleated surfaces $\Sigma_n$ in $M$ which {\em exits} the end $e_1$: Each $\Sigma_n$ is contained in 
 an isolating neighborhood $E_1$ of $e_1$ 
 and for every compact subset $K\subset M$ there exists $n_0$ such that $\Sigma_n\cap K=\emptyset$ for all $n\ge n_0$.

 \medskip
 {\bf Case 1:} The image $\Ga_1$ of $\pi_1(S_1)\to \pi_1(\bar{M})=\Ga$ has infinite index (equivalently, index $\ge 3$) in $\Ga$. 
 In other words, the covering map
 $$
 \H^3/\Ga_1\to M
 $$
 has infinite multiplicity. Then, by the Thurston--Canary covering theorem, \cite{C2},  the manifold $M_1=\H^3/\Ga_1$ 
 has at least one geometrically finite end (more precisely, all ends not covering $e_1$ are geometrically finite). In other words, $\La(\Ga_1)\ne S^2$. 
 By  Theorem \ref{BJnon}, the Hausdorff dimension of the non-conical limit set of $\Ga_1$ equals $2$. 
 Hence, $\dim(\Lambda_{nc}(\Gamma))\geq \dim(\Lambda(\Gamma_{1}))=2$ and, thus, $\dim(\Lambda_{nc}(\Gamma))=2$. 

\medskip 
 {\bf Case 2:} After passing to an index two subgroup, we may (and will) assume that $\Ga_1=\Ga$. Then $\bar{M}$ is a compression body, obtained by 
 attaching 2-handles and 3-handles to $S_1\times [0,1]$ along $S_1\times \{0\}$. The assumption that $\Ga$ is not free means that $\bar{M}$ is not a handlebody, i.e. $\bar{M}$ has at least one more boundary component besides $S_1$; equivalently, $M$ has at least two ends. 
If at least one of the ends of $M$ is geometrically finite then $\La(\Ga_1)\ne S^2$, and the same argument as in Case 1 concludes the proof. 
Hence, we will assume that all ends of $M$ are geometrically infinite.  As we discussed in section \ref{sec:ends-3D}, 
 there exists another complete hyperbolic  manifold $M'=\H^3/\Ga'$ without cusps, homeomorphic to  $M$ such that:
 
 a. All ends of $M'$ besides the end $e'_1$ corresponding to $e_1$, are geometrically finite. 
 
 b. There is a diffeomorphism $h: M'\to M$ which is bi-Lipschitz on a neighborhood of the end $e'_1$; more precisely, 
 there exist closed isolating neighborhoods  $E_1, E_1'$ of $e_1, e_1'$ respectively, such that the restriction $h: E'_1\to E_1$ is 
 bi-Lipschitz with respect to the path-metrics of $E_1, E_1'$ induced from $M, M'$ respectively.

The assumption that the injectivity radius of $M$ at $E_1$ is bounded from below and (b) imply that the injectivity radius of $M'$ is bounded from 
below at $E'_1$. Since all other ends of $M'$ are geometrically finite, it follows that the injectivity radius of $M'$ is bounded from below. 
In particular, according to \cite{BJ}, $\dim \La_{nc}(\Ga')=2$. We will use this and the bi-Lipschitz diffeomorphism 
$h: E'_1\to E_1$ to conclude that $\dim \La_{nc}(\Ga)>0$.

Since the injectivity radius of $M'$ is bounded from below, 
the neighborhood $E'_1$ is $a$-narrow for some $a$. 

Since $\Ga_1=\Ga$, $E_1, E'_1$ have connected preimages $\tilde E_1, \tilde E'_1$ in $\H^3$  (under the universal covering maps $\H^3\to M, \H^3\to M'$). 
 We  fix a lift $\tilde h$ of the bi-Lipschitz map $h: E'_1\to E_1$. 
For a base-point $x'\in \tilde E'_1$, we will be using the base-point $x= \tilde h(x')$ in $\tilde E_1$. 

After enlarging $E'_1$ (and using the fact that all the ends of $M'$ besides $e'_1$ are geometrically finite), 
we may assume that the neighborhood $E'_1$ of $e_1'$ is the convex core $M'_c$ of $M'$. Thus, $\tilde E'_1=C(\La(\Ga'))$, the closed convex hull of the limit set of $\Ga'$. In particular, the intrinsic path-metric of $\tilde E'_1$ equals the hyperbolic distance function restricted from $\H^3$.

For $\tilde{E}_1$ this is, of course, false, but every its intrinsic geodesic disjoint from 
$\partial \tilde E_1$  is a geodesic in $\H^3$ as well. In view of the quasiisometry $\tilde h: \tilde E'_1 \to \tilde E_1$, there exists $\delta\in \R$ such that 
the path-metric space  $\tilde E_1$ is $\delta$-hyperbolic, i.e. every (intrinsic) geodesic triangle in $\tilde E_1$ is $\delta$-slim.

Consider an (intrinsic) geodesic triangle $[abc]$ in $\tilde{E}_1$. If all three sides of this triangle are disjoint from $\partial \tilde E_1$ then the intrinsic Gromov-product $(a|b)_c$ equals the extrinsic Gromov-product defined via the metric on $\H^3$. Since $[abc]$ is $\delta$-slim, it suffices to assume that the geodesics $ac, bc$ are disjoint from the (intrinsic) $\delta$-neighborhood of $\partial \tilde E_1$. In view of Lemma \ref{lem:visual}, we then obtain:

\begin{lemma}\label{lem:bh} 
If 
$\xi_1, \xi_2\in \La(\Ga)$ are $(x,\delta)$-deep limit points 
with respect to $E_1$, then the visual angle (computed from $x$) between $\xi_1, \xi_2$ is uniformly bi-H\"older to the Gromov-distance $d_x(\xi_1, \xi_2)$ computed with respect to the intrinsic metric on $\tilde E_1$. 
Here the H\"older constants depend only on $x$ and not on $\xi_1, \xi_2$.    
\end{lemma} 

\begin{lemma}\label{lem:holder} 
There exists a constant $D_1$ such that for any $D> D_1$ the map $\tilde h: \tilde E_1'\to \tilde E_1$ extends to a bi-H\"older embedding 
$$
\theta: \La_{(x',D)}(E'_1)\to \La(e_1). 
$$
\end{lemma}
\proof For each geodesic ray $\rho'=x'\xi'$ in $\tilde E'_1$, the composition $\tilde h\circ \rho'$ is an intrinsic 
$(L,0)$-quasigeodesic in $\tilde E_1$ where $L$ is the  bi-Lipschitz constant of $\tilde h: \tilde E_1'\to \tilde E_1$. 
Therefore, by the Morse Lemma, the image of this quasigeodesic is within distance $D_0:=D(L,0,\delta)$ from an intrinsic geodesic $\rho$ in $\tilde E_1$ 
with the same origin $x=\tilde{h}(x')$.

For any $D\ge 0$, if the minimal distance between the image of $\rho'$ and $\partial \tilde E'_1$ is 
$\ge D$ then the minimal distance between the image of $\rho$ and $\partial E_1$ is $\ge L^{-1}D - D_0$. In particular, for $D> L(\delta + D_0)$,  
whenever $x'\xi'$ is disjoint from the closed $D$-neighborhood of $\partial \tilde E'_1$, the geodesic $\rho$ is disjoint from the closed $\delta$-neighborhood of 
$\partial \tilde E_1$.  In other words, the intrinsic geodesic $\rho$ is also a geodesic $x\xi$ in $\H^3$ which avoids the intrinsic 
closed $\delta$-neighborhood of $\partial \tilde E_1$. Furthermore, if $\rho'$ projects to a proper ray in $E'_1$, the quasiray 
 $\tilde h\circ \rho'$ projects to a proper quasiray in $E_1$ and so does the ray $\rho$. By Lemma \ref{lem:limit-of-narrow-end}, 
 the point $\xi$ is necessarily a limit point of $\Ga$, i.e. belongs to $\La(e_1)$. We thus obtain a map
$$
\theta: \La_{(x',D)}(E'_1)\to \La(e_1), \quad \xi'\mapsto \xi.  
$$

We now prove a bi-H\"older estimate where $D$ is chosen so that $D> D_1=L(\delta + D_0)$. 
 First of all, according to Theorem \ref{holder}, the quasiisometry $\tilde h: \tilde E_1'\to \tilde E_1$ of intrinsic metrics induces a bi-H\"older bijection of the Gromov-boundaries. By convexity, the identity map  $\tilde E_1'\to \tilde E_1'$ is a isometry from the extrinsic metric to the intrinsic metric, hence, it induces a bi-H\"older homeomorphism between the Gromov boundaries equipped, respectively, with the angular metric and a visual metric defined via Gromov-product. By Lemma \ref{lem:bh}, for pairs of intrinsic geodesic rays $(r_1, r_2)$ in $\tilde E_1$ emanating from $x$ and 
which are disjoint from the closed $\delta$-neighborhood of $\partial \tilde E_1$, the angle between $r_1, r_2$ at $x$ is uniformly bi-H\"older to the visual distance between the corresponding points of the Gromov-boundary of $\tilde E_1$ (equipped with the intrinsic metric). We conclude that the map $\theta$ 
is bi-H\"older to its image. \qed 

\medskip 
We now can conclude the proof of Theorem \ref{dim2}. 
Recall that $E'_1$ is $a$-narrow. Hence, by Corollary \ref{cor:edim}, for $x'\in \tilde E'$ satisfying 
$$
d(x', \partial \tilde E'_1)>  a+\delta_0
$$
and $D= d(x', \partial \tilde E'_1)-  a-\delta_0$,
$$
\dim \La(e'_1)= \dim \La_{(x',D)}(E'_1). 
$$
We choose $x'$ so that 
$$
d(x', \partial \tilde E'_1)>  D_1+ a+\delta_0,   
$$ 
where $D_1$ is the constant given by Lemma \ref{lem:holder}. 
Then, by Theorem \ref{BJnon}, $2=\dim \La(e'_1)=\dim \La_{(x',D)}(E'_1)$, while Lemma \ref{lem:holder}  implies that $\theta(\La_{(x',D)}(E'_1))\subset \La(e_1)$ has positive Hausdorff dimension.  \qed

\end{document}